\documentclass{article}
\usepackage{
amsmath,
amsfonts,
latexsym, 
amsrefs,
amssymb
}

\usepackage{bm}

\usepackage{tocloft}

\usepackage{cite}
\usepackage{subcaption}

\usepackage{pgfplots}
\usepackage{tikz}
\usepackage{tikz,fullpage}
\usetikzlibrary{arrows,%
                petri,%
                topaths, automata}%
\usepackage{tkz-berge} 

\usepackage{graphicx}
\newcommand{\bla}{\bm{\lambda}}

\DeclareMathOperator{\Prob}{\mathbb{P}}   

\usepackage{enumerate}
\usepackage{mathrsfs}
\usepackage{wrapfig}

\usepackage{color}

\numberwithin{equation}{section}
\newcommand{\boeta}{\bm{\eta}}

\newcommand{\rd}{{\rm d}}

\newcommand{\bu}{{\bf{u}}}

\newcommand{\bw}{{\bf{w}}}

\newcommand{\al}{\alpha}

\newcommand{\be}{\begin{equation}}
\newcommand{\ee}{\end{equation}}

\newcommand{\e}{{\varepsilon}}

\newcommand{\la}{\lambda}

\usepackage{amsmath} 
\usepackage{amssymb}
\usepackage{amsthm}

\renewcommand{\b}[1]{\bm{\mathrm{#1}}} 

\newcommand{\wt}{\widetilde}

\newcommand{\ii}{\mathrm{i}} 

\renewcommand{\epsilon}{\varepsilon}
\renewcommand{\leq}{\leqslant}
\renewcommand{\geq}{\geqslant}

\renewcommand{\le}{\leq}
\renewcommand{\ge}{\geq}

\renewcommand{\P}{\mathbb{P}}
\newcommand{\E}{\mathbb{E}}

\DeclareMathOperator{\OO}{O}
\DeclareMathOperator{\oo}{o}

\theoremstyle{plain} 
\newtheorem{theorem}{Theorem}[section]
\newtheorem*{theorem*}{Theorem}
\newtheorem{lemma}[theorem]{Lemma}
\newtheorem*{lemma*}{Lemma}

\newtheorem*{corollary*}{Corollary}

\newtheorem*{proposition*}{Proposition}

\newtheorem{definition}[theorem]{Definition}
\newtheorem*{definition*}{Definition}

\newtheorem*{example*}{Example}

\newtheorem*{remark*}{Remark}
\newtheorem*{remarks*}{Remarks}

\makeatletter
\renewcommand{\subsection}{\@startsection
{subsection}
{2}
{0mm}
{-\baselineskip}
{0 \baselineskip}
{\normalfont\bf\itshape}} 
\makeatother

\usepackage{dsfont}
\usepackage{stmaryrd}

\newcommand{\nc}{\normalcolor}

\newcommand{\bq}{{\bf q}}

\renewcommand{\b}[1]{\boldsymbol{\mathrm{#1}}} 

\def\@empty{}

\def\author#1{\par
    {\centering{\authorfont#1}\par\vspace*{0.05in}}
}

\def\titlefont{\fontsize{13}{15}\bfseries\boldmath\selectfont\centering{}}
\def\authorfont{\fontsize{13}{15}}
\def\abstractfont{\fontsize{8}{10}}

\let\affiliationfont\rhfont

\def\address#1{\par
    {\centering{\affiliationfont#1\par}}\par\vspace*{11pt}
}

\def\keywords#1{\par
    \vspace*{8pt}
    {\authorfont{\leftskip18pt\rightskip\leftskip
    \noindent{\it\small{Keywords}}\/:\ #1\par}}\vskip-12pt}

\def\title#1{
    \thispagestyle{plain}
    \vspace*{-14pt}
    \vskip 79pt
    {\centering{\titlefont #1\par}}%
    \vskip 1em
}

\renewenvironment{abstract}{\par%
    \vspace*{6pt}\noindent 
    \abstractfont
    \noindent\leftskip18pt\rightskip18pt
}{%
  \par}

\renewcommand{\b}[1]{\boldsymbol{\mathrm{#1}}} 

\usepackage{tocloft}

\makeatletter
\renewcommand{\section}{\@startsection
{section}
{1}
{0mm}
{-2\baselineskip}
{1\baselineskip}
{\normalfont\large\scshape\centering}} 
\makeatother

\begin{document}

~\vspace{-1.5cm}

\title{Random band matrices}

\vspace{0.7cm}
\noindent\hspace{4cm}\begin{minipage}[b]{0.5\textwidth}

 \author{P. Bourgade}

\address{Courant Institute, New York University\\
  bourgade@cims.nyu.edu}
 \end{minipage}
\begin{minipage}[b]{0.5\textwidth}

 \end{minipage}


\begin{abstract}
We survey recent mathematical results about the spectrum of random band matrices.
We start by exposing the Erd{\H o}s-Schlein-Yau dynamic approach, its application to Wigner matrices, and extension to other mean-field models. We then introduce random band matrices and the problem of their Anderson transition.
We finally describe a method to obtain delocalization and universality in some sparse  regimes, highlighting the role of quantum unique ergodicity.
\end{abstract}


{\let\thefootnote\relax\footnotetext{\noindent This work is supported by the NSF grant DMS\#1513587.}}

\tableofcontents

\vspace{0.5cm}

This review explains the interplay between eigenvectors and eigenvalues statistics in random matrix theory, when the considered models are not of mean-field type, meaning that the interaction is short range and  geometric constraints enter in the definition of the model.

If the range or strength of the interaction is small enough, it is expected that eigenvalues statistics will fall into the {\it Poisson universality class}, intimately related to the notion of independence. Another class emerged in 
the past fifty years for many correlated systems, initially from calculations on random linear operators.
This {\it random matrix universality class} was proposed by  Wigner \cite{Wig1957}, first as a model for stable energy levels of typical heavy nuclei.
The models he introduced have since been understood to connect to integrable systems, growth models, analytic number theory and multivariate statistics (see e.g. \cite{Dei2017}).

Ongoing efforts to understand universality classes  are essentially of two types. First,
{\it integrability} consists in finding possibly new statistics for a few models, with methods including combinatorics and representation theory. Second,
{\it universality} means enlarging the range of models with random matrix statistics, through probabilistic methods. For example, the Gaussian random matrix ensembles are mean-field integrable models, from which local spectral statistics  can be established for the more general Wigner matrices,  by comparison, as explained in Section 1. For random operators with shorter range, no integrable models are known, presenting a major difficulty in understanding whether their spectral statistics will fall in the Poisson or random matrix class.

In Wigner's original theory, the eigenvectors play no role. However, their statistics are essential in view of a famous dichotomy of spectral behaviors, widely studied since Anderson's tight binding model \cite{And}: 
\begin{enumerate}[(i)]
\item Poisson spectral statistics usually occur together with localized eigenstates,
\item
random matrix eigenvalue distributions should coincide with delocalization of eigenstates.
\end{enumerate}
The existence of the localized phase has been established for the Anderson model in any dimension \cite{FroSpe}, but delocalization has remained elusive for all operators relevant in physics. An important question consists in proving extended states and GOE local statistics for one such model\footnote{GOE eigenvalues statistics appear  in Trotter's tridiagonal model \cite{Tro}, which is clearly local, but the entries need varying variance adjusted to a specific profile.}, giving theoretical evidence for conduction in solids.
How localization implies Poisson statistics is well understood, at least for the Anderson model \cite{Min}. In this note, we explain the proof of a strong notion of delocalization (quantum unique ergodicity), and how it implies random matrix spectral statistics, for the $1d$ random band matrix (RBM) model.

This model can be defined for general dimension  ($d=1,2,3$): vertices are elements of $\Lambda=\llbracket 1,N\rrbracket^d$
and $H=(H_{ij})_{i,j\in\Lambda}$ have centered real entries, independent 
up to the symmetry $H_{ij}=H_{ji}$. The band width $W<N/2$ means 
\begin{equation}
H_{ij}=0\ \mbox{ if}\  |i-j|>W,\label{eqn:defband}
\end{equation}
where $|\cdot|$ is the periodic ${\rm L}^1$ distance on $\Lambda$,
and all non-trivial $H_{ij}$'s have a variance $\sigma_{ij}^2$ with the same order of magnitude, normalized by $\sum_j \sigma_{ij}^2=1$
for any $i\in\Lambda$. Mean-field models correspond to $W=N/2$. When $W\to\infty$, the empirical spectral measure of $H$ converges to the semicircle distribution $\rd \rho_{\rm sc}(x)=\frac{1}{2\pi}(4-x^2)^{1/2}\rd x$.

It has been conjectured that the random band matrix model exhibits the localization-delocalization (and Poisson-GOE)  transition at some critical band width $W_c(N)$ for eigenvalues in the bulk of the spectrum $|E|<2-\kappa$.  The localized  regime supposedly occurs for $W\ll W_c$ and delocalization for $W\gg W_c$, where
\begin{equation}\label{eqn:transition exponents}
W_c=
\left\{
\begin{array}{ll}
N^{1/2}&\mbox{for}\,\, d=1,\\
(\log N)^{1/2}&\mbox{for}\,\, d=2,\\
\OO(1)&\mbox{for}\,\, d=3.
\end{array}
\right.
\end{equation}
This transition corresponds to localization length $\ell\approx W^2$ in dimension 1, $\ell\approx e^{W^2}$ in dimension 2.

\begin{figure}[h]
\centering
~\vspace{1cm}
\includegraphics[width=4cm]{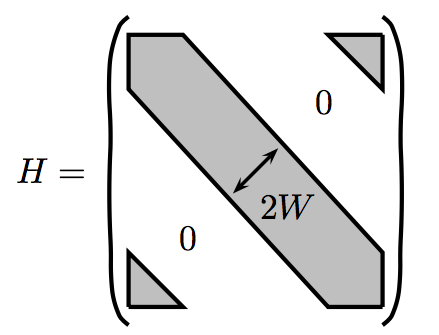}
\hspace{0.5cm}
\begin{subfigure}{.3\textwidth}
~\vspace{2.5cm}
\begin{tikzpicture}[scale=0.6]
\begin{axis}[axis lines=middle, enlargelimits=true, xtick=\empty, ytick=\empty,xmin=0,xmax=4,ymin=0,ymax=1,height=6cm,width=10cm]
\addplot [dashed,domain=0:4,samples=100,line width=2pt] {exp(-x)};
\addplot [domain=0:4, samples=100,line width=2pt] {2*x*exp(-x^2)};
\end{axis}
\end{tikzpicture}
\end{subfigure}
\hspace{0.1cm}
\begin{subfigure}{.3\textwidth}
~\vspace{2.5cm}
\begin{tikzpicture}[scale=0.6]
\begin{axis}[axis lines=middle, enlargelimits=true, xtick=\empty,label style={font=\Large\bf}, ytick=\empty,xmin=0,xmax=4,ymin=0,ymax=1,height=6cm,width=10cm,xlabel={$\alpha$},ylabel={$({\rm Av} |u_k|^2)(\alpha)$},xlabel style={at=(current axis.right of origin), anchor=west,right=3.3cm},ylabel style={at=(current axis.above origin), anchor=south,above=1.7cm}]
\addplot [dashed,domain=0:4, samples=100,line width=2pt] {0.9*exp(-20*((x-2.55)^4-0.1*(x-2.5)^2)};
\addplot [domain=0:4, samples=100,line width=2pt] {0.2};
\end{axis}
\node at (1,-0.1) {{{\small$1$}}};
\node at (7.5,-0.1) {{{\small$N$}}};
\node at (0.1,1) {{{\small$1/N$}}};
\node at (5.4,-0.3) {{{\small$\ell\approx W^2$}}};
\draw[<->,line width=0.5pt](3.7,0.1) -- (6.7,0.1);
\end{tikzpicture}
\end{subfigure}

\vspace{-3.5cm}
\caption{Conjectural behavior of the RBM model for $d=1$. For any eigenvalue $|\lambda_k|<2-\kappa$, the rescaled  gap $N\rho_{\rm sc}(\la_k)(\la_{k+1}-\la_k)$  converges to an exponential random variable for $W\ll N^{1/2}$, and the Gaudin GOE distribution for $W\gg N^{1/2}$. The associated eigenvector $u_k$ is localized on $\ell\approx W^2$ sites for  $W\ll N^{1/2}$, it is  flat for $W\gg N^{1/2}$. Here $({\rm Av} f)(\alpha)=(2n)^{-1}\sum_{|i-\alpha|<n}f(i)$ where $1\ll n\ll W^2$ is some averaging scale.}\label{fig:d=1}
\end{figure}

\vspace{-0.05cm}

This review first explains universality techniques for mean-field models. We then state recent progress for random band matrices, including the existence of the delocalized phase for $d=1$ \cite{BouErdYauYin2017,BouYauYin2018}, explaining how quantum unique ergodicity is proved by dynamics. We finally explain, at the heuristic level, a connection between quantum unique ergodicity for band matrices and the Gaussian free field, our main goal being to convince the reader that the transition exponents in (\ref{eqn:transition exponents}) are natural.\\

For the sake of conciseness, we only consider the orthogonal symmetry class corresponding to  random symmetric matrices with real entries. Analogous results hold in the complex Hermitian class.

\section{Mean-field random matrices}

\subsection{Integrable model.}\ 
The Gaussian orthogonal ensemble (GOE) consists in the probability density 
\begin{equation}\label{eqn:GOE}
\frac{1}{Z_N}e^{-\frac{N}{4}{\rm Tr}(H^2)}
\end{equation}
with respect to the Lebesgue measure on the set on $N\times N$ symmetric matrices. This corresponds to all entries being Gaussian and independent up to the symmetry condition, with off-diagonal entries $H_{ij}\sim N^{-1/2}\mathscr{N}(0,1)$, and diagonal entries $H_{ii}\sim(N/2)^{-1/2}\mathscr{N}(0,1)$.

Our normalization is chosen so that the eigenvalues $\la_1\leq \dots\leq \la_N$ (with associated eigenvectors $u_1,\dots,u_N$) have a converging empirical measure: $\frac{1}{N}\sum_{k=1}^N\delta_{\la_i}\to \rd\rho_{{\rm sc}}$ almost surely. A more detailed description of the spectrum holds at the microscopic scale,  in the bulk and at the edge: there exists a translation invariant point process $\chi_1$ \cite{MehGau} and a distribution ${\rm TW}_1$ (for Tracy and Widom \cite{TraWid}) such that 
\begin{align}
\label{sine1}\sum_{k=1}^N\delta_{N\rho_{\rm sc}(E)(\la_k-E)}&\to \chi_1,\\
\label{TW1}N^{2/3}(\la_N-2)&\to {\rm TW}_1,
\end{align}
in distribution. Note that $\chi_1$ is independent of $E\in(-2+\kappa,2-\kappa)$, for any fixed, small, $\kappa>0$.

Concerning the eigenvectors, for any $\mathscr{O}\in O(N)$, from (\ref{eqn:GOE})  the distributions of $\mathscr{O}^{\rm t}H\mathscr{O}$ and $H$ are the same, so that the eigenbasis $\bu=(u_1,\dots,u_N)$ of $H$ is Haar-distributed (modulo a sign choice) on $O(N)$: $\mathscr{O}\bu$ has same distribution as $\bu$. In particular, any $u_k$ is uniform of the sphere $\mathscr{S}^{(N-1)}$, and has the same distribution as $\mathscr{N}/\|\mathscr{N}\|_2$ where $\mathscr{N}$ is a centered Gaussian vector with covariance ${\rm Id}_N$. This implies that for any deterministic sequences of indices $k_N\in \llbracket 1,N\rrbracket$ and  unit vectors $\bq_N\in \mathscr{S}^{(N-1)}$ (abbreviated $k,\bq$), the limiting {\it Borel-L\'evy law} holds:
\begin{equation}\label{eqn:Gaussian}
N^{1/2}\langle u_k,\bq\rangle\to \mathscr{N}(0,1)
\end{equation}
in distribution.
This microscopic behavior can be extended to several projections being jointly Gaussian.

The  fact that eigenvectors are extended can be  quantified in different manners. For example, for the GOE model, for any small $\e>0$ and large $D>0$, we have
\begin{equation}\label{eqn:deloc}
\mathbb{P}\left(\|u_k\|_\infty\geq \frac{N^\e}{\sqrt{N}}\right)\leq N^{-D},
\end{equation}
which we refer to as {\it delocalization}  (the above $N^\e$ can also be replaced by some logarithmic power). 

Delocalization does not imply that the eigenvectors are flat in the sense of Figure \ref{fig:d=1}, as $u_k$ could be supported on a small fraction of $\llbracket 1,N\rrbracket$. 
A strong notion of flat eigenstates was introduced by Rudnick and Sarnak \cite{RudSar1994} for Riemannian manifolds: they conjectured that for any negatively curved and compact $\mathcal{M}$ with volume measure $\mu$,
$$
\int_A|\psi_k(x)|^2\mu(\rd x)\underset{k\to\infty}{\longrightarrow}\int_A\mu(\rd x),
$$
for any $A\subset \mathcal{M}$.
Here $\psi_k$ is an eigenfunction (associated to the eigenvalue $\lambda_k$) of the Laplace-Beltrami operator, $0\leq \lambda_1\leq\dots\leq \lambda_k\leq\dots$ and $\|\psi_k\|_{{\rm L}^2(\mu)}=1$.  
This {\it quantum unique ergodicity} (QUE) notion 
strengthens the quantum ergodicity proved in \cite{Shn1974,Col1985,Zel1987}, defined by an additional averaging on $k$ and proved for a wide class of manifolds and deterministic regular graphs  \cite{AnaLeM2013} (see also \cite{BroLin}).
QUE was rigorously proved for arithmetic surfaces,
\cites{Lin2006,Hol2010,HolSou2010}. 
We will consider a probabilistic version of QUE at a local scale, for eigenvalues in the bulk of the spectrum. By simple properties of the uniform measure on the unit sphere
it is clear that the following version holds for the GOE:
for any given (small) $\e>0$ and (large) $D>0$,
for $N\geq N_0(\e,D)$,
for any deterministic sequences $k_N\in \llbracket \kappa N,(1-\kappa)N\rrbracket$ and $I_N\subset\llbracket 1,N\rrbracket$ (abbreviated $k,I$),   we have
\be\label{eqn:localque}
\Prob\left(\left|\sum_{\alpha\in I} (u_k(\alpha)^2-\frac{1}{N})\right|>\frac{N^\e|I|^{1/2}}{N}\right)\leq N^{-D}.
\ee
We now consider the properties (\ref{sine1}), (\ref{TW1}) (\ref{eqn:Gaussian}), (\ref{eqn:deloc}), (\ref{eqn:localque}) for the following general model.

\begin{definition}[Generalized Wigner matrices]\label{generalizedWigner}
A sequence $H_N$ (abbreviated $H$) of real symmetric centered random matrices is a generalized Wigner matrix if there exists $C,c>0$ such that $\sigma_{ij}^2:={\rm Var}(H_{ij})$ satisfies
\begin{equation}\label{eqn:meanfield}
c\leq N \sigma_{ij}^2\leq C\ \mbox{for all $N,i,j$ and}\ \sum_j\sigma_{ij}^2=1\ \mbox{for all $i$.}
\end{equation}
We also assume  subgaussian decay of the distribution of $\sqrt{N}H_{ij}$, uniformly in $i,j,N$, for convenience (this could be replaced by a finite high moment assumption).
\end{definition}

\subsection{Eigenvalues universality.}\ The second constraint in (\ref{eqn:meanfield}) imposes the macroscopic behavior of the limiting spectral measure: $\frac{1}{N}\sum_{k=1}^N\delta_{\la_i}\to \rd\rho_{{\rm sc}}$ for all generalized Wigner matrices. This convergence to the semicircle distribution was strengthened up to optimal polynomial scale, thanks to an advanced diagrammatic analysis of the resolvent of $H$.

\begin{theorem}[Rigidity of the spectrum \cite{ErdYauYin2012Rig}]\label{thm:rigidity}
Let $H$ be a generalized Wigner matrix as in Definition \ref{generalizedWigner}. Define $\hat k=\min(k,N+1-k)$ and $\gamma_k$ implicitly by $\int_{-2}^{\gamma_k}\rd\rho_{\rm sc}=\frac{k}{N}$. Then for any $\e>0$, $D>0$ there exists $N_0$ such that for $N>N_0$, $k\in\llbracket 1,N\rrbracket$, we have
\begin{equation}\label{eqn:rig}
\mathbb{P}\left(|\la_k-\gamma_k|>N^{-\frac{2}{3}+\e}(\hat k)^{-\frac{1}{3}}\right)\leq N^{-D}.
\end{equation}
\end{theorem}

Given the above scale of fluctuations, a natural problem consists in the limiting distribution.
In particular, the (Wigner-Dyson-Mehta) conjecture states that (\ref{sine1})
holds for random matrices way beyond the integrable GOE class. It has been proved in a series of works in the past
years, with important new techniques based on the Harish-Chandra-Itzykson-Zuber integral \cite{Joh2001} (in the special case of Hermitian symmetry class), 
the dynamic interpolation through Dyson Brownian motion \cite{ErdSchYau2011II} and the Lindeberg exchange principle \cite{TaoVu2011}. The initial universality statements for general classes required an averaging over the energy level $E$ \cite{ErdSchYau2011II} or the first four moments of the matrix entries to match the Gaussian ones \cite{TaoVu2011}.

We aim at explaining the dynamic method which was applied in a remarkable variety of settings. For example, GOE local eigenvalues statistics hold for generalized Wigner matrices.

\begin{theorem}[Fixed energy universality \cite{BouErdYauYin2014}]\label{thm:fixedenergy} The convergence 
 (\ref{sine1}) holds for generalized Wigner matrices.
\end{theorem}


\noindent The key idea for the proof, from \cite{ErdSchYau2011II}, is interpolation through matrix Dyson Brownian motion (or its Ornstein Uhlenbeck version) 
\begin{equation}\label{eqn:matrixDBM}
\rd H_t=\frac{1}{\sqrt{N}}\rd B_t-\frac{1}{2}H_t\rd t
\end{equation}
with initial condition $H_0=H$, where $(B_{ij})_{i<j}$ and $(B_{ii}/\sqrt{2})_i$ are independent standard Brownian motions. The GOE measure (\ref{eqn:GOE}) is the equilibrium for these dynamics. 
The proof proceeds in two steps, in which the dynamics
$(H_t)_{t\geq 0}$
is analyzed through complementary viewpoints. One relies on the repulsive eigenvalues dynamics, the other on the matrix structure. Both steps require some a priori knowledge on eigenvalues density, such as Theorem \ref{thm:rigidity}.\\

\noindent {\it First step: relaxation.} For any $t\geq N^{-1+\e}$,  (\ref{sine1}) holds:
$\sum_{k=1}^N\delta_{N\rho_{\rm sc}(E)(\la_k(t)-E)}\to \chi_1$, where 
we denote $\la_1(t)\leq\dots\leq\la_N(t)$ the eigenvalues of $H_t$. The proof relies on the
Dyson Brownian motion for the eigenvalues dynamics \cite{Dys}, given by
\begin{equation}\label{eqn:DBM}
\rd\la_k(t)=\frac{\rd \wt B_{k}(t)}{\sqrt{N}}+\left(\frac{1}{N}\sum_{\ell\neq k}\frac{1}{\la_k(t)-\la_\ell(t)}-\frac{1}{2}\la_k(t)\right)\rd t
\end{equation}
where the $\wt B_k/\sqrt{2}$'s are standard Brownian motions. Consider the dynamics (\ref{eqn:DBM}) with a different
initial condition $x_1(0)\leq \dots\leq x_N(0)$ given by the eigenvalues of a GOE matrix. By taking the difference between these two coupled stochastic differential equations we observe that $\delta_\ell(t):=e^{t/2}(x_\ell(t)-\la_\ell(t))$ satisfy an integral equation of parabolic type \cite{BouErdYauYin2014}, namely
\begin{equation}\label{aftercoupling}
\partial_t \delta_\ell(t) = 
\sum_{k\neq \ell} b_{k\ell}(t)(\delta_k(t)-\delta_\ell(t)),\ \ b_{k\ell}(t)=\frac{1}{N(x_\ell(t)-x_k(t))(\la_\ell(t)-\la_k(t))}.
\end{equation}
From Theorem \ref{thm:rigidity}, in the bulk of the spectrum we expect that $b_{k\ell}(t)\approx N/(k-\ell)^2$, so that 
H\"{o}lder regularity holds for $t\gg N^{-1}$: $\delta_k(t)=\delta_{k+1}(t)(1+{\rm o}(1))$, meaning $\la_{k+1}(t)-\la_k(t)
=y_{k+1}(t)-y_{k}(t)+{\rm o}(N^{-1})$. Gaps between the $\la_k$'s and $x_k$'s therefore become identical, hence equal to the GOE gaps as the law of 
$y_{k+1}(t)-y_{k}(t)$ is invariant in time. In fact, an equation similar to (\ref{aftercoupling}) previously appeared 
in the first proof of GOE gap statistics for generalized Wigner matrices \cite{ErdYau2012singlegap}, emerging from a Helffer-Sj\"{o}strand representation instead of a probabilistic coupling.
Theorem \ref{thm:fixedenergy} requires a much more precise analysis of (\ref{aftercoupling}) \cite{BouErdYauYin2014,LanSosYau2016}, but the conceptual picture is clear from the above probabilistic coupling of eigenvalues.

Relaxation after a short time can also be understood by functional inequalities for relative entropy \cite{ErdSchYau2011II,ErdYauYin2012Univ}, a robust method which also gives  GOE statistics when averaging over the energy level $E$. In the special case of the Hermitian symmetry class, relaxation also follows from explicit formulas for the eigenvalues density at time $t$ \cite{Joh2001,ErdPecRamSchYau2010,TaoVu2011}.\\

\noindent{\it Second step: density}. For any $t\leq N^{-\frac{1}{2}-\e}$, $\sum_{k=1}^N\delta_{N\rho_{\rm sc}(E)(\la_k(t)-E)}$ and $\sum_{k=1}^N\delta_{N\rho_{\rm sc}(E)(\la_k(0)-E)}$ have the same distribution at leading order.
This step can be proved by a simple It{\^o} lemma based on the matrix evolution \cite{BouYau2017}, which takes a particularly simple form for Wigner matrices (i.e. $\sigma_{ij}^2=N^{-1}+N^{-1}\mathds{1}_{i=j}$). It essentially states that
for any smooth function $F(H)$
we have
\begin{equation}\label{eqn:constant}
 \E F\left(H_t\right)  -  \E F\left(H_0\right)
 = \OO( t {N^{1/2}}) \sup_{i\leq j,0\leq s\leq t} \E \left((N^{3/2} |H_{ij}(s)^3| + \sqrt N  |H_{ij}(s)|) \big | \partial_{{ij}}^3 F  (H_s  )\big |\right)
\end{equation}
where $\partial_{ij}=\partial_{H_{ij}}$.
In particular, if $F$ is stable in the sense that $\partial_{{ij}}^3 F =\OO(N^\e)$ with high probability (this is known for functions encoding the microscopic behavior thanks to the a-priori rigidity estimates from Theorem \ref{thm:rigidity}), 
then the same local statistics as for $t=0$ holds up to time $N^{-\frac{1}{2}-\e}$.

Invariance of local spectral statistics has also been proved by other methods, for example by a reverse heat flow  when the entries have a smooth enough density \cite{ErdSchYau2011II}, or the Lindeberg exchange principle \cite{TaoVu2011} for matrices with moments of the entries coinciding up to fourth moment.

\subsection{Eigenvectors universality.}\ 
Eigenvalues rigidity (\ref{eqn:rig}) was an important estimate for the proof of Theorem \ref{thm:fixedenergy}. Similarly, to understand the eigenvectors distribution, one needs to first identify their natural fluctuation scale. By analysis of the resolvent of $H$, the following was first proved when $\bq$ is an element from the canonical basis \cite{ErdSchYau2009,ErdYauYin2012Rig}, and extended to any direction.

\begin{theorem}[Isotropic delocalization \cite{KnoYin2013II,BloErdKnoYauYin2014}]
For any sequence of generalized Wigner matrices,  $\e,D>0$, there exists $N_0(\e,D)$ such that for any $N\geq N_0$, deterministic  $k$ and  unit vector $\bq$, we have
$$
\mathbb{P}\left(\langle u_k,\bq\rangle\geq N^{-\frac{1}{2}+\e}\right)\leq N^{-D}.
$$\end{theorem}
The more precise fluctuations (\ref{eqn:Gaussian}) were proved by the Lindeberg exchange principle in \cite{KnoYin2013,TaoVu2012}, under the assumption of the first four (resp. two) moments of $H$ matching the Gaussian ones, for eigenvectors associated to the spectral bulk (resp. edge). This  L\'evy-Borel law holds without these moment matching assumptions, and some form of quantum unique ergodicity comes with it.

\begin{theorem}[Eigenvectors universality and weak QUE \cite{BouYau2017}]\label{thm:BY}
For any sequence of generalized Wigner matrices,  
and any deterministic  $k$ and  unit vector $\bq$, the convergence (\ref{eqn:Gaussian}) is true.

Moreover, for any $\e>0$ there exists $D>0$ such that
 (\ref{eqn:localque}) holds.
\end{theorem}

The above statement is a {\it weak} form of QUE, holding for some small $D=D(\e)$ although
it should be true for any large $D>0$. Section 3 will show a strong form of QUE
for some band matrices.

The proof of Theorem \ref{thm:BY} follows the dynamic idea already described for eigenvalues, by considering the evolution of the eigenvectors through (\ref{eqn:matrixDBM}). The {\it density} step is similar: with (\ref{eqn:constant}) one can 
show that the distribution of $\sqrt{N}\langle u_k(t),\bq\rangle$ is almost invariant up to time $t\leq N^{-\frac{1}{2}-\e}$.
The {\it relaxation} step is significantly different from the coupling argument described previously. The eigenvectors dynamics are given by
$$
\rd u_k=\frac{1}{\sqrt{N}}\sum_{\ell\neq k}\frac{\rd \wt B_{k\ell}}{\lambda_k-\lambda_\ell}u_\ell
-\frac{1}{2N}\sum_{\ell\neq k}\frac{\rd t}{(\la_k-\la_\ell)^2}u_k,
$$
where the $\wt B_{k\ell}$'s are independent standard Brownian motions, and most importantly  independent from 
the $\wt B_k$'s from (\ref{eqn:DBM}).
This eigenvector flow was computed   
in the context of Brownian motion on ellipsoids \cite{NorRogWil1986},   real Wishart processes \cite{Bru1989}, and for GOE/GUE in \cite{AndGuiZei2010}.

Due to its complicated structure and high dimension, this eigenvector flow had not been previously analyzed. 
Surprisingly, these  dynamics can be reduced to a 
multi-particle random walk in a dynamic random environment. 
More precisely, let a configuration $\boeta$ consist in $d$ points of $\llbracket1,N\rrbracket$, with possible repetition. The number of particles at site $x$ is $\eta_x$.
A configuration obtained by moving a particle from $i$ to $j$ is denoted $\boeta^{ij}$. The main observation from \cite{BouYau2017} is as follows.
First denote $z_k=\sqrt{N}\langle \bq, u_k\rangle$, which is random and time dependent. 
Then associate to a configuration $\boeta$ with $j_k$ points at $i_k$, the renormalized moments observables (the $\mathscr{N}_{i_k}$ are independent Gaussians)
conditionally to the eigenvalues path,
\begin{equation}\label{eqn:obs1}
f_{t,\bla}(\boeta)=
\E\left(\prod_{k=1}^N z_{i_k}^{2j_k} \mid\bla\right)/\,\E\left(\prod_{k=1}^N\mathscr{N}_{i_k}^{2j_k}\right).
\end{equation}
Then $f_{t,\lambda}$ satisfies the parabolic partial differential equation

\begin{minipage}[r]{0.4\linewidth}
\includegraphics[width=7cm]{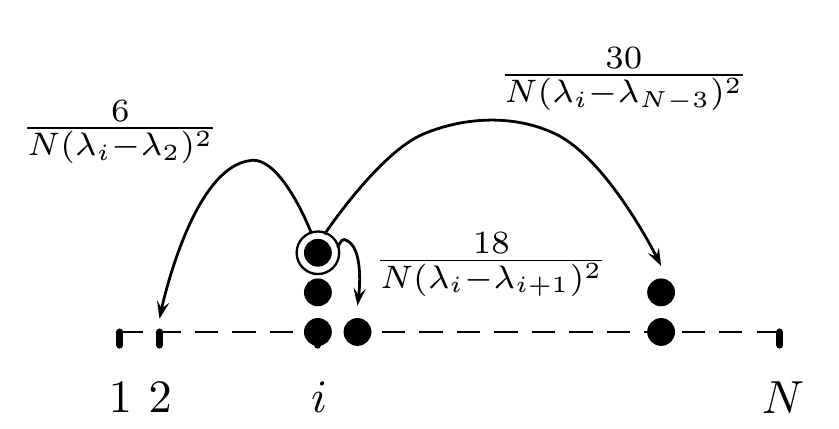}
\end{minipage}
\hspace{1cm}
\begin{minipage}{0.5\linewidth}
\begin{equation}\label{eqn:EMF}
\hspace{-5cm}\partial_t f_{t,\bla}(\boeta)=\mathscr{B}(t)  f_{t,\bla}(\boeta) 
\end{equation}
where
$$
\hspace{-1cm}\mathscr{B}(t)  f(\boeta) =\frac{1}{N} \sum_{i \neq j} 2 \eta_i (1+ 2 \eta_j) \frac{f(\boeta^{ij})-f(\boeta)}{(\lambda_i(t)-\lambda_j(t))^2}.
$$
\end{minipage}

As shown in the above drawing, the generator $\mathscr{B}(t)$ corresponds to a random walk on the space of configurations $\boeta$, with time-dependent rates given by the eigenvalues dynamics.
This equation is parabolic and by the scale argument explained for (\ref{aftercoupling}), $f_{t,\bla}$ becomes locally constant (in fact, equal to 1 by normalization constraint) for $t\geq N^{-1+\e}$.
This H\"{o}lder regularity is  proved by a maximum principle.

\subsection{Other models.}\ 
The  described dynamic approach applies beyond generalized Wigner matrices.
We do not attempt to give a complete list of applications of this method. Below are a few results.

\begin{enumerate}[(i)]

\item
Wigner-type matrices designate variations of Wigner matrices with non centered $H_{ii}$'s  \cite{LeeSchSteYau2015}, or the normalization constraint in (\ref{eqn:meanfield}) not satisfied (the limiting spectral measure differs from semicircular) \cite{AjaErdKru2015}, or the $H_{ij}$'s non-centered and correlated \cite{AjaErdKru2018,ErdKruSch2017,Che}. In all cases, GOE bulk statistics are known.

\item
Random graphs also have bulk or edge GOE statistics when the connectivity grows fast enough with $N$, as proved for example for the Erd{\H o}s-Renyi \cite{ErdKnoYauYinER,LanHuaYau2015,HuaLanYau2017,LeeSch2015} and uniform $d$-regular models  \cite{BauHuaKnoYau2017}. Eigenvectors statistics are also known to coincide with the GOE for such graphs \cite{BouHuaYau2017}.

\item
For $\beta$-ensembles,  the external potential does not impact local statistics, a fact first shown when $\beta=1,2,4$ (the classical invariant ensembles)
by asymptotics of orthogonal polynomials \cite{BleIts1999,  Dei1999, DeiGio,Lub2009,PasShc1997}. The dynamics approach extended this result to any $\beta$ \cite{BouErdYau2014,BouErdYau2011}.
Other methods based  on sparse models \cite{KriRidVir} and transport maps  \cite{BekFigGui2013,Shc2013} were also applied to either $\beta$-ensembles or multimatrix models \cite{FigGuiII}.

\item
The convolution model  $D_1+U^*D_2U$, where $D_1$, $D_2$ are diagonal and $U$ is uniform on ${\rm O}(N)$, appears in free probability theory. Its empirical spectral measure in understood up to the optimal scale  \cite{BaoErdSch2017}, and GOE bulk statistics were proved in \cite{CheLan2017}.

\item For small mean-field perturbations of diagonal matrices (the Rosenzweig-Porter model), GOE statistics \cite{LanSosYau2016} occur with localization \cite{Ben2017,VonWar2017}. We refer to \cite{FacVivBir} for the physical meaning of this unusual regime.

\item
Extremal statistics. The smallest gaps in the spectrum of Gaussian ensembles and  Wigner matrices have the same law \cite{Bou2018}, when the matrix entries are smooth. The relaxation step (\ref{aftercoupling}) was quantified with an optimal error  so that the smallest spacing scale ($N^{-4/3}$ in the GUE case \cite{BenBou2011}) can  be perceived.
\end{enumerate}

The above matrix models are mean-field, a constraint inherent to the  proof strategy previously described. Indeed, the density step requires the matrix entries to fluctuate substantially: lemmas of type (\ref{eqn:constant}) need a constant variance of the entries along the dynamics (\ref{eqn:matrixDBM}).

\section{Random band matrices and the Anderson transition}

In the Wigner random matrix model, the entries, which represent the quantum transition rates between two quantum states, are all of comparable size.
More realistic models involve  geometric structure, as typical quantum transitions only occur between nearby states.
In this section we briefly review key results for  Anderson and band matrix models.

\subsection{Brief and partial history of random Schr\"{o}dinger operators.}\ 
Anderson's random Schr\"{o}dinger operator  \cite{And} on $\mathbb{Z}^d$ describes a  system  
 with spatial structure. It is of type 
 \begin{equation}\label{eqn:RS}
H_{\rm RS}=\Delta+\lambda V
\end{equation}
 where $\Delta$ is the discrete Laplacian and  the random variables $V(x)$, $x\in\mathbb{Z}^d$, are i.i.d and centered with variance $1$. 
 The parameter $\lambda>0$ measures the strength of the disorder.
 The spectrum of $H_{\rm RS}$ is supported on $[-2d,2d]+\lambda{\rm supp}(\mu)$ where $\mu$ is the distribution of $V(0)$

 Amongst the many mathematical contributions to this model, Anderson's initial motivation (localization, hence the suppression of electron transport due to disorder) was proved rigorously by 
   Fr\"{o}hlich and Spencer \cite{FroSpe} by a multiscale analysis: localization holds
 for strong disorder or at energies
where the density of states $\rho(E)$ is small (localization for a related one-dimensional model was previously proved by Golsheid, Molchanov and Pastur \cite{Gol}).
An alternative derivation was given in Aizenman and Molchanov \cite{AizMol}, who introduced a fractional moment method. 
From the scaling theory of localization \cite{AbrAndLicRam},  extended states supposedly occur in dimensions $d\geq 3$ for $\lambda$ small enough, while eigenstates are only marginally localized for $d=2$.

Unfortunately, there has been no progress in establishing the delocalized regime for the random Schr\"odinger operator on $\mathbb{Z}^d$. 
The existence of absolutely continuous spectrum (related to extended states) in the presence of substantial disorder  is only known when $\mathbb{Z}^d$ is replaced by homogeneous trees  \cite{Kle1994}.

These results and conjecture were initially for the Anderson model in infinite volume. If we denote $H_{\rm RS}^N$ the operator (\ref{eqn:RS}) restricted to the box $\llbracket-N/2,N/2\rrbracket^d$ with periodic boundary conditions, its spectrum still lies on a compact set and one expects that the bulk eigenvalues in the microscopic scaling (i.e. multiplied by $N^d$) converge to either Poisson or GOE statistics ($H_{\rm RS}^N$
corresponds to GOE rather than GUE because it is a real symmetric matrix).
Minami proved Poisson spectral statistics  from exponential decay of the resolvent \cite{Min}, in cases where localization  in infinite volume is known.
For $H_{\rm RS}^N$, not only is the existence of delocalized states in dimension three open, but also there is no clear understanding about how extended states imply GOE spectral statistics.

\subsection{Random band matrices: analogies, conjectures, heuristics.}\ 
The band matrix model we will consider was essentially already defined around (\ref{eqn:defband}). In addition, 
in the following we will assume  subgaussian decay of the distribution of $W^\frac{d}{2}H_{ij}$, uniformly in $i,j,N$, for convenience (this could be replaced by a finite high moment assumption).

Although random band matrices and the random Schr\"{o}dinger operator (\ref{eqn:RS}) 
are different, they are both local (their matrix elements $H_{ij}$ vanish
when $|i-j|$ is large). The models are expected to have the same 
properties when
\begin{equation}\label{eqn:analogy}
\lambda\approx \frac{1}{W}.
\end{equation}
For example, eigenvectors for the Anderson model in one dimension are proved to decay
exponentially fast with a localization length proportional to $\lambda^{-2}$, in agreement with the analogy (\ref{eqn:analogy}) and the conjecture (\ref{eqn:transition exponents}) when $d=1$.
 For $d=2$, it is conjectured that all states are localized with a
localization length of order $\exp(W^2)$ for band matrices,  $\exp(\lambda^{-2})$ for the Anderson model, again coherently with (\ref{eqn:analogy}) and (\ref{eqn:transition exponents}).
For some mathematical justification of the analogy (\ref{eqn:analogy}) from the point of view of perturbation theory, we refer to \cite[Appendix 4.11]{Spe2012}.

The origins of conjecture (\ref{eqn:transition exponents})  first lie on  numerical evidence, at least for $d=1$.
In \cite{ConJ-Ref1} it was observed, based on computer simulations, that the bulk eigenvalue statistics and eigenvector localization length of $1d$ random band matrices are essentially a function of $W^2/N$,
with the sharp transition as in (\ref{eqn:transition exponents}).  Fyodorov and Mirlin gave the first theoretical explanation for this transition \cite{FyoMir}. They considered a slightly different ensemble with complex Gaussian entries decaying exponentially fast at distance greater than $W$ from the diagonal. Based on a non-rigorous supersymmetric approach \cite{Efe1997}, they approximate relevant random matrix statistics with expectations for a related  $\sigma$-model, from which a saddle point method gives the localization/delocalization transition for $W\approx \sqrt{N}$. Their work also gives an estimate on the localization length $\ell$, anywhere in the spectrum \cite[equation (19)]{FyoMir}: it is expected that at energy level $E$ (remember our normalization $\sum_j\sigma_{ij}^2=1$ for $H$ so that the equilibrium measure is $\rho_{\rm sc}$),
$$
\ell\approx \min(W^2(4-E^2),N).
$$
With this method, they were also able to conjecture explicit formulas for the distribution of
eigenfunction components and related quantities for any scaling ratio $W^2/N$ \cite{FyoMir1994}.

Finally, heuristics for localization/delocalization transition exponents follow from the conductance fluctuations theory developed by Thouless \cite{Thouless}, based on scaling arguments. For a discussion of mathematical aspects of the Thouless criterion, see \cites{Spe2012,Spe}, and 
 \cite[Section III]{Wang} for some rigorous scaling theory of localization.  This  criterion was introduced in the context of Anderson localization, and was  applied  in \cite{Sod2010,Sod2014} to $1d$ band matrices, including at the edge of the spectrum, in agreement with the prediction from \cite{FyoMir}. A different heuristic argument for (\ref{eqn:transition exponents}) is given in Section 3, for any dimension in the bulk of the spectrum.

\subsection{Results.}\ 
The density of states ($\E\left(N^{-1}\sum_k\delta_{\lambda_k}\right)$) of properly scaled random band matrices in dimension $1$ converges to the semicircular distribution for any $W\to\infty$, as proved in \cite{BogMolPas1991}.
This convergence was then strengthened and fluctuations around the semicircular law were studied in \cite{Gui,AndZei2006,JanSahSos2016,LiSos} by the method of moments, at the macroscopic scale.

  Interesting transitions extending the microscopic one (\ref{eqn:transition exponents}) are supposed to occur at mesoscopic scales $\eta$, giving a full phase diagram in $(\eta,W)$. The work  \cite{ErdKno2015I} rigorously analyzed parts of this diagram by studying linear statistics in some mesoscopic range and in any dimension, also by a moment-based approach.
  
The miscroscopic scale transitions (\ref{eqn:transition exponents}) are harder to understand, but recent progress  allowed to prove the existence of localization and delocalization for some polynomial scales in $W$.
These results are essentially of four different types: $(i)$ the localization side for general models, $(ii)$ localization and delocalization for specific Gaussian models, $(iii)$ delocalization for general models.
Finally, $(iv)$ the edge statistics are fully understood by the method of moments. Unless otherwise stated, all results below are restricted to $d=1$.\\

\noindent {\it (i) Localization for general models.}  A seminal result in the analysis of random band matrices is the following estimate on the localization scale. For simplicity one can assume that the entries of $H$ are i.i.d. Gaussian, but the method from \cite{Sch} allows to treat more general distributions.

\begin{theorem}[The localization regime for band matrices \cite{Sch}]
Let $\mu>8$. There exists $\tau>0$ such that for large enough $N$, for any $\alpha,\beta\in\llbracket 1,N\rrbracket$ one has
$$
\mathbb{E}\left(\sup_{1\leq k\leq N}|u_k(\alpha)u_k(\beta)|\right)\leq W^\tau e^{-\frac{|\alpha-\beta|}{W^\mu}}.
$$
\end{theorem}

\noindent Localization therefore holds simultaneously for all eigenvectors when $W\ll N^{1/8}$, which was improved to $W\ll N^{1/7}$ in \cite{PelSchShaSod} for some specific Gaussian model  described below.\\

\noindent {\it (ii) Gaussian models with specific variance profile and supersymmetry.}  
For some Gaussian band matrices, the supersymmetry (SUSY) technique
gives a purely analytic derivation of spectral properties.
This approach  
has first been developed by physicists  \cite{Efe1997}. 
A rigorous  supersymmetry method started with the expected density of states on arbitrarily short scales for a $3d$ band matrix ensemble 
\cite{DisPinSpe2002}, extended to $2d$ in \cite{DisLag2017} (see \cite{Spe2012} for much more about the mathematical aspects of SUSY).  More recently, the work  \cite{Shc} proved local 
  GUE local statistics for $W\geq c N$, and  delocalization was obtained in a strong sense for individual
  eigenvectors,  when 
  $W\gg N^{6/7}$ and the first four moments of the matrix entries match the Gaussian ones \cite{BaoErd2015}.  
These recent rigorous results assume complex entries and hold for $|E|<\sqrt{2}$, for a block-band structure of the matrix with a specific variance profile.

We briefly illustrate the SUSY method for moments of the characteristic polynomial: remarkably, this is currently the only observable for which the transition at $W\approx\sqrt{N}$ was proved.
Consider a  matrix $H$ whose entries are complex centered Gaussian variables  such that 
$$
\mathbb{E}(H_{ij}H_{\ell k})=\mathds{1}_{i=k,j=\ell}J_{ij}\ \mbox{where}\ J_{ij}=(-W^2\Delta+1)^{-1}_{ij},
$$
and $\Delta$ is the discrete Laplacian on $\llbracket 1,N\rrbracket$ with periodic boundary condition. The variance $J_{ij}$ is exponentially small for $|i-j|>W^{1+\e}$, so that $H$ can be considered a random band matrix with band width $W$.
Define 
$$
F_2(E_1,E_2)=\mathbb{E}\left(\det(E_1-H)\det(E_2-H)\right),\ D_2=F_2(E,E).
$$
\begin{theorem}[Transition for characteristic polynomials \cite{Sch1,SchMT}]\label{thmS}
For any $E\in(-2,2)$ and $\e>0$, we have
$$
\lim_{N\to\infty} (D_2)^{-1}F_2\left(E+\frac{x}{N\rho_{\rm sc}(E)},E-\frac{x}{N\rho_{\rm sc}(E)}\right)=
\left\{
\begin{array}{ll}
1&\mbox{if}\ N^\e<W<N^{\frac{1}{2}-\e}\\
\frac{\sin(2\pi x)}{2\pi x}&\mbox{if}\ N^{\frac{1}{2}+\e}<W<N
\end{array}
\right..
$$
\end{theorem}
Unfortunately, currently the local eigenvalues statistics cannot be identified from products of characteristic polynomials: they require ratios which are more difficult to analyze by the SUSY method.

We briefly mention the key steps of the proof of Theorem \ref{thmS}. First, an integral representation for $F_2$ is obtained by integration over Grassmann variables. 
These variables give convenient formulas for the product of characteristic polynomials: they allow to express the determinant as a Gaussian-type integral. Integrate over the Grassmann variables then gives an integral representation (in complex variables) of the moments of interest. More precisely, the  Gaussian representation for $F_2\left(E+\frac{x}{N\rho_{\rm sc}(E)},E-\frac{x}{N\rho_{\rm sc}(E)}\right)$, from \cite{Sch1}, is
$$
\frac{1}{(2\pi)^{N}}\frac{1}{\det{J}^2}\int e^{-\frac{W^2}{2}\sum_{j=-n+1}^n{\rm Tr}(X_j-X_{j-1})^2
-\frac{1}{2}\sum_{j=-n}^n{\rm Tr}(X_j+\frac{\ii \Lambda_E}{2}+\ii \frac{\Lambda_x}{N\rho_{\rm sc}(E})^2}\prod_{j=-n}^n
\det(X_j-\ii \Delta_E/2){\rm d}X_j,
$$
where $N=2n+1$, $\Delta_E={\rm diag}(E,E)$, $\Delta_x={\rm diag}(x,-x)$, and $\rd X_j$ is the Lebesgue measure on $2\times 2$ Hermitian matrices. This form of the correlation of characteristic polynomial is then analyzed by steepest descent.
Analogues of the above representation hold in any dimension, where the matrices $X_j,X_k$, are coupled in a quadratic way when $k$ and $j$ are neighbors in $\mathbb{Z}^d$, similarly to the Gaussian free field.

Finally, based on their integral representations, it is expected that random band matrices behave like  $\sigma$-models, which are used by physicists to understand complicated statistical mechanics systems. We refer to the recent work \cite{Sch2018} for rigorous results in this direction.\\

\noindent{\it (iii)\ Delocalization for general models.}
Back to general models with no specific distribution of the entries (except sufficient decay of the distribution, for example subgaussian), the first delocalization results for random band matrices relied on a difficult 
analysis of their resolvent.

For example, the Green's function   was controlled
 down to the scale $W^{-1}$ in \cite{ErdYauYin2012Univ}, implying that
 the localization length of all eigenvectors is at least $W$. 
 Analysis of the resolvent also gives full delocalization  for most eigenvectors, for $W$ large enough. In the theorem below, 
we say that an eigenvector $u_k$ is subexponentially localized at scale $\ell$ if there exists $\e>0$, $I\subset\llbracket 1,N\rrbracket$, $|I|\leq \ell$, such that $\sum_{\alpha\not\in I}|u_k(\alpha)|^2<e^{-N^\e}$.

\begin{theorem}[Delocalized regime on average \cite{HeMarc2018}]\label{thm:averaged}
Assume $W\gg N^{7/9}$ and $\ell\ll N$. Then the fraction of eigenvectors subexponentially localized on scale $\ell$ vanishes as $N\to\infty$, with large probability.
\end{theorem}

\noindent This result for $W\geq N^{6/7}$ was previously obtained in \cite{ErdKno2013}, for $W\geq N^{4/5}$  in \cite{ErdKnoYauYin}, and similar statements were proved in higher dimension.

Delocalization was recently proved without averaging, together with eigenvalues statistics and flatness of individual eigenvectors. 
The main new ingredient is that quantum unique ergodicity is a convenient delocalization notion, proved by dynamics.

To simplify the statement below, assume that $H$ is a Gaussian-divisible ,
in the sense that for $|i-j|\leq W$, $\sqrt{W}H_{ij}$ is the sum of two independent random variables, $X+\mathscr{N}(0,c)$, where $c$ is an arbitrary small constant
(the result holds for more general entries).

\begin{theorem}[Delocalized regime \cite{BouYauYin2018}]\label{thm:notaveraged}
Assume $W\gg N^{3/4+a}$ for some $a>0$. Let  $\kappa>0$ be fixed.
\begin{enumerate}[(a)]
\item For any $E\in(-2+\kappa,2-\kappa)$ the eigenvalues statistics at energy level $E$ converge to the GOE, as in  (\ref{sine1}).
\item The bulk eigenvectors are  delocalized: for any (small) $\e>0$, (large) $D>0$, for $N\geq N_0(\e,D,\kappa)$ and  $k\in\llbracket \kappa N,(1-\kappa)N\rrbracket$, we have 
$$
\mathbb{P}\left(\|u_k\|_\infty>N^{-\frac{1}{2}+\e}\right)<N^{-D}.
$$
\item The bulk eigenvectors are flat on any scale greater than $W$. More precisely,
for any given (small) $\e>0$ and (large) $D>0$,
for $N\geq N_0(\e,D,\kappa)$,
for any deterministic $k\in \llbracket \kappa N,(1-\kappa)N\rrbracket$ and interval $I\subset\llbracket 1,N\rrbracket$, $|I|>W$,  we have
$$
\Prob\left(\left|\sum_{\alpha\in I} (u_k(\alpha)^2-\frac{1}{N})\right|>N^{-\frac{3}{2}a+\e} \frac{|I|}{N}\right)\leq N^{-D}.
$$
\end{enumerate}

\end{theorem}

A strong form of QUE similar to $(c)$ holds for random $d$-regular graphs \cite{BauHuaYau2017}, the proof relying on exchangeability. For models with geometric constraints, other ideas are explained in the next section.

Theorem \ref{thm:notaveraged} relies on a mean-field reduction strategy initiated in \cite{BouErdYauYin2017}, and an extension of the dynamics (\ref{eqn:EMF}) to observables much more general than (\ref{eqn:obs1}), as explained in the next section. New ingredients compared to Theorem \ref{thm:averaged} are  (a) quantum unique ergodicity for mean-field models after Gaussian perturbation, in a strong sense, (b) estimates on the resolvent of the band matrix at the (almost macroscopic) scale $N^{-\e}$. 

The current main limitation of the method to approach the transition $W_c=N^{1/2}$ comes from (b). These resolvent estimates are obtained by intricate diagrammatics developed in a series of previous works including  \cite{ErdKnoYauYin}, extended to generalized resolvents and currently only proved for $W\gg N^{3/4}$ \cite{BouFanYauYin2018,FanYin2018}.\\

\noindent{\it (iv)\ Edge statistics.}\  The transition in eigenvalues statistics is understood at the edge of the spectrum: the prediction from the Thouless criterion was made rigorous by a subtle method of moments. 
This was proved under the assumption that $\sqrt{2W}(H_{ij})_{i\leq j}$ are $\pm1$ independent centered Bernoulli random variables, but the method applies to more general distributions.

\begin{theorem}[Transition at the edge of the spectrum \cite{Sod2010}] 
If $W\gg N^{\frac{5}{6}}$, then (\ref{TW1}) holds. If $W\ll N^{\frac{5}{6}}$, (\ref{TW1}) does not hold.
\end{theorem}

\noindent Finally, for eigenvectors (including at the edge of the spectrum), localization cannot hold on less than $W/\log(N)$ entries as proved in \cite{BenPec2014}, also by the method of moments.

\section{Quantum unique ergodicity and universality}

For non mean-field models, eigenvalues and eigenvectors interplay extensively, and their statistics should be understood jointly. 
Localization (decay of Green's function) is a useful a priori estimate in the proof of Poisson statistics for the Anderson model \cite{Min}, and in a similar way we explain below why quantum unique ergodicity implies GOE statistics.

\subsection{Mean-field reduction.}\ 
The method  introduced in \cite{BouErdYauYin2017} for GOE statistics of band matrices proceeds as follows.
We 
decompose  the $1d$ band matrix from (\ref{eqn:defband}) and its eigenvectors as
$$
   H=  \begin{pmatrix} A  & B^* \cr B & D  \end{pmatrix}, \quad 
   \bu_j:=   \begin{pmatrix}\b w_j \cr \b p_j \end{pmatrix},
$$
where $A$ is a $W\times   W$ matrix. From the eigenvector equation $H \bu_j = \lambda_j \bu_j$ we have
$
   (A- B^* \frac{1}{D - \lambda _j }B) \bw_j  = \lambda_j  \bw_j.
$
The matrix
elements of  $A$ do not vanish and thus the above eigenvalue  problem  features a mean-field random matrix 
(of smaller size). 
Hence one can considers the eigenvector equation 
$
Q_e \bw_k(e) = \xi_k(e)\bw_k(e) 
$
where 
\begin{equation}\label{eqn:complement}
Q_e= A-B^* (D-e)^{-1}B, 
\end{equation}
and  $\lambda_k(e)$, $\bw_k (e)$ $(1\leq k\leq W)$ are  eigenvalues and normalized eigenvectors.
As illustrated below, the slopes of the functions $e\mapsto\lambda_k(e)$ seem to be locally equal and concentrated:
$$
\frac{\rd}{\rd e}\lambda_k(e)\approx1-\frac{1}{\sum_{\alpha=1}^W w_k(\alpha)^2}(1+\oo(1))\approx 1-\frac{N}{W},
$$ 
which holds for $e$ close to $\lambda_k$. The first equality is a simple perturbation formula\footnote{The perturbation formula gives a slightly different equation, replacing $\bw_k$ by the eigenvector of a small perturbation of $H$, but we omit this technicality.}, and the second  is true provided QUE for $\bu_k$ holds, in the sense of equation (\ref{eqn:localque}) for example.

\begin{figure}
\centering
\begin{subfigure}{.4\textwidth}
  \centering
\begin{tikzpicture}
\node[anchor=south west,inner sep=0] (x) at (0,0) {\includegraphics[width=7cm]{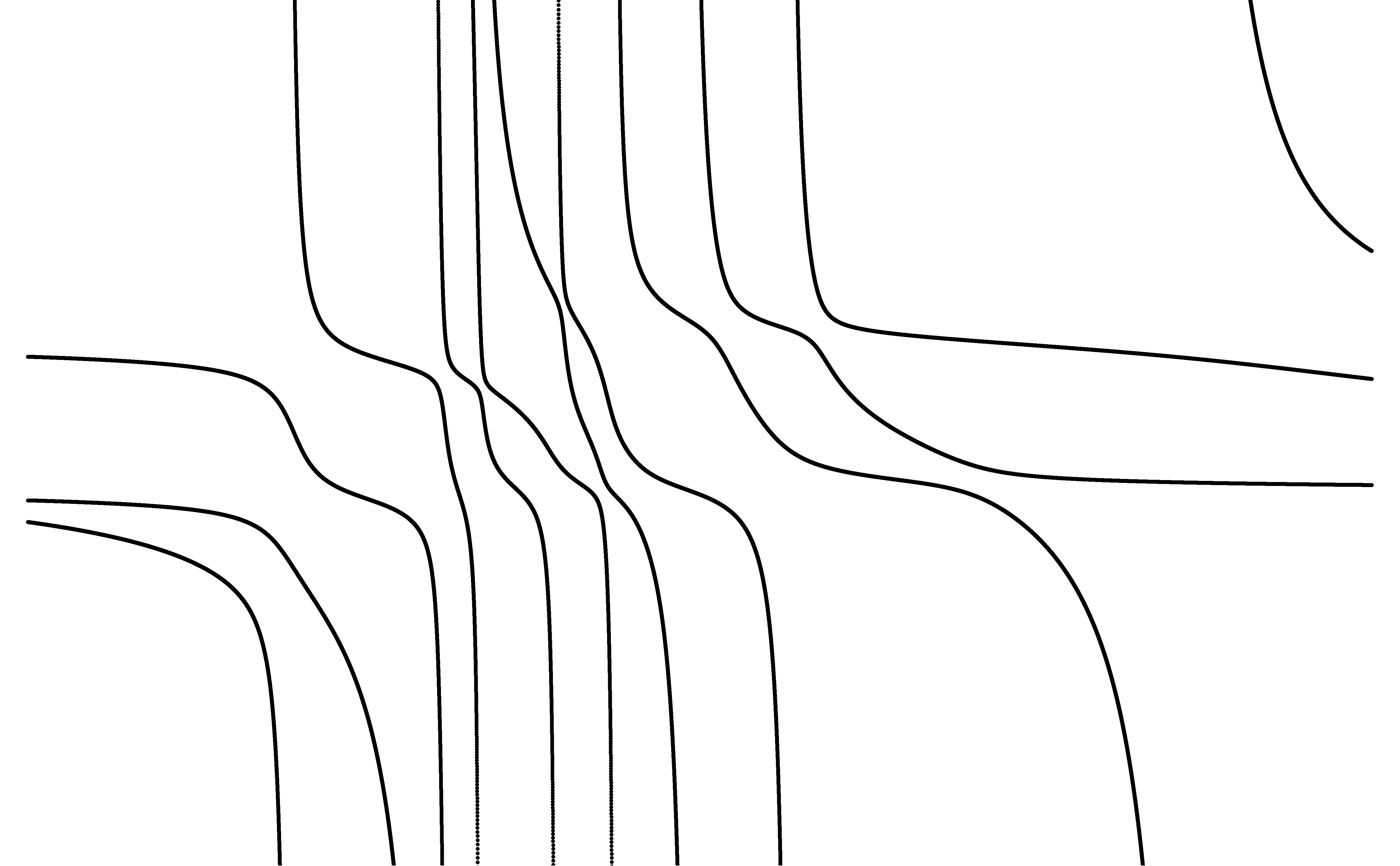}};
\begin{scope}[x={(x.south east)},y={(x.north west)}]
\draw[black,thick,rounded corners] (0.5,0.55) rectangle (0.64,0.75);
\draw [black,->] (0,0.5) -- (1.02,0.5);
\draw [black,->] (0.47,0) -- (0.47,1);
\draw [black,dashed] (0.47,0.5) -- (0.81,1);
\draw [black,dashed] (0.47,0.5) -- (0.13,0);
\fill[black] (0.198,0.1)  circle[black,radius=1.5pt];
\fill[black] (0.26,0.19)  circle[black,radius=1.5pt];
\fill[black] (0.313,0.265)  circle[black,radius=1.5pt];
\fill[black] (0.339,0.305)  circle[black,radius=1.5pt];
\fill[black] (0.388,0.378)  circle[black,radius=1.5pt];
\fill[black] (0.423,0.426)  circle[black,radius=1.5pt];
\fill[black] (0.435,0.445)  circle[black,radius=1.5pt];
\fill[black] (0.452,0.468)  circle[black,radius=1.5pt];
\fill[black] (0.52,0.58)  circle[black,radius=1.5pt];
\fill[black] (0.555,0.622)  circle[black,radius=1.5pt];
\fill[black] (0.585,0.665)  circle[black,radius=1.5pt];
\fill[black] (0.198,0.1)  circle[black,radius=1pt];
\node[black] at (1.05,0.5) {$e$};
\end{scope}
\end{tikzpicture}
  \caption{A simulation of  eigenvalues of $Q_e=A-B^*(D-e)^{-1}B$, i.e. functions $e\mapsto \la_j(e)$. Here $N=12$ and $W=3$. The $\lambda_i$'s are the abscissa of the intersections with the diagonal.}
\end{subfigure}%
\begin{subfigure}{.1\textwidth}
  \centering
\begin{tikzpicture}
\draw [black,dashed] (0,0) -- (0,0);
\end{tikzpicture}
\end{subfigure}%
\begin{subfigure}{.4\textwidth}
\centering
\vspace{-1.1cm} \hspace{-2cm}
\begin{tikzpicture}
\node[anchor=south west,inner sep=0] (y) at (0,0) {\includegraphics[width=8.5cm]{C.jpg}};
\begin{scope}[x={(y.south east)},y={(y.north west)}]
\draw[white,ultra thick,fill=white]  (0,0) rectangle (1,1);
\draw [black,-,thick,dashed] (0.21,0.16) -- (0.65,0.8);
\draw [black,-,thick,dashed] (0.37,0.16) -- (0.37,0.8);
\draw [black,-,thick] (0.1,0.7) -- (0.8,0.3);
\draw [black,-,thick] (0.1,0.8) -- (0.8,0.4);
\draw [black,-,thick] (0.1,0.93) -- (0.8,0.53);
\draw [black,-,thick] (0.1,1) -- (0.8,0.6);
\draw [black,-,thick] (0.1,0.64) -- (0.8,0.24);
\draw [black,-,thick] (0.1,0.5) -- (0.8,0.1);
\draw [black,-,thick] (0.1,0.4) -- (0.8,0.02);
\draw [black,-,thick] (0.1,0.32) -- (0.63,0.04);
\fill[black] (0.26,0.235)  circle[black,radius=1.8pt];
\fill[black] (0.3,0.29)  circle[black,radius=1.8pt];
\fill[black] (0.35,0.36)  circle[black,radius=1.8pt];
\fill[black] (0.415,0.46)  circle[black,radius=1.8pt];
\fill[black] (0.445,0.5)  circle[black,radius=1.8pt];
\fill[black] (0.495,0.573)  circle[black,radius=1.8pt];
\fill[black] (0.56,0.665)  circle[black,radius=1.8pt];
\fill[black] (0.595,0.715)  circle[black,radius=1.8pt];
\draw[black,fill=black] (0.37,0.775)  +(-1.3pt,-1.3pt) rectangle +(1.3pt,1.3pt) ;
\draw[black,fill=black] (0.37,0.645)  +(-1.3pt,-1.3pt) rectangle +(1.3pt,1.3pt) ;
\draw[black,fill=black] (0.37,0.545)  +(-1.3pt,-1.3pt) rectangle +(1.3pt,1.3pt) ;
\draw[black,fill=black] (0.37,0.485)  +(-1.3pt,-1.3pt) rectangle +(1.3pt,1.3pt) ;
\draw[black,fill=black] (0.37,0.345)  +(-1.3pt,-1.3pt) rectangle +(1.3pt,1.3pt) ;
\draw[black,fill=black] (0.37,0.255)  +(-1.3pt,-1.3pt) rectangle +(1.3pt,1.3pt) ;
\draw[black,fill=black] (0.37,0.18)  +(-1.3pt,-1.3pt) rectangle +(1.3pt,1.3pt) ;
\draw [black,->](0.37,0.845)-- (0.495,0.774);
\filldraw[white,fill=white]
(0,0) -- (0,1) -- (1,1) -- (1,0) -- cycle
(0.2,0.15) -- (0.65,0.15) -- (0.65,0.8) -- (0.2,0.8) -- cycle;
\draw[black,thick,rounded corners] (0.2,0.15) rectangle (0.65,0.8);
\draw [black,->,thick] (0.37,0.05) -- (0.37,0.13);
\node[black] at (0.37,0.02) {$e$};

 \draw [black,->,thick] (0.31,0.05) -- (0.31,0.13);
\node[black] at (0.31,0.02) {$\lambda'$};

\draw [black,->,thick] (0.26,0.05) -- (0.26,0.13);
\node[black] at (0.26,0.02) {$\lambda$};

\draw [black,->](0.37,0.775)-- (0.475,0.715);
\draw [black,->] (0.37,0.645) -- (0.44,0.605);
\draw [black,->] (0.37,0.545)-- (0.415,0.52);
\draw [black,->] (0.37,0.49)-- (0.398,0.468);
\draw [black,->] (0.37,0.255) -- (0.33,0.274);
\draw [black,->] (0.37,0.18) -- (0.31,0.209);
\end{scope}
\end{tikzpicture}
  \caption{Zoom into the framed region of Figure (a), for large $N,W$: the curves $\la_j(e)$ are almost parallel, with  slope about $1-N/W$.  The eigenvalues of  $A-B^*(D-e)^{-1}B$ and those of $H$ are related by a
  projection to the diagonal followed by a projection to
  the horizontal axis.
 }
\end{subfigure}

\caption{The idea of mean-field reduction: universality of gaps between eigenvalues for fixed $e$ implies universality on the diagonal through parallel projection.
For $e$ fixed, we label the curves by $\la_k(e)$.  }
\label{Fig1}
\end{figure}

The GOE local spectral statistics hold for $Q_e$ in the sense (\ref{sine1}) (it is a mean-field matrix so results from \cite{LanSosYau2016} apply), hence it also holds for $H$ by parallel projection: GOE local spectral statistics follow from QUE.

This reduces the problem to QUE for band matrices, which is proved by the same mean-field reduction strategy: on the one hand, by choosing different overlapping blocks $A$ along the diagonal, QUE for $H$
follows from QUE for $Q_e$ by a simple patching procedure (see section 3.3 for more details); on the other hand, QUE for mean-field models is known thanks to a strengthening of the eigenvector moment flow method \cite{BouYau2017,BouHuaYau2017}, explained below.

\subsection{The eigenvector moment flow.}\ In this paragraph, $(u_k)_k$ now refers to the eigenvectors of a $N\times N$ mean-field random matrix, with eigenvalues $(\la_k)_k$, as in Section 1.

Obtaining quantum unique ergodicity from the regularity of equation (\ref{eqn:EMF}) (the eigenvector moment flow) is easy: $\sqrt{N}\langle\bq,u_k\rangle$ has limiting Gaussian moments for any $\bq$, hence the entries of $u_k$ are asymptotically independent Gaussian and the following variant of (\ref{eqn:localque}) holds for $Q_e$ by Markov's inequality ($w_k$ is rescaled to a unit vector):
there exists $\e>0$ such that
for any  deterministic $1\leq k\leq N$ and $I\subset\llbracket 1,N\rrbracket$, for any $\delta>0$ we have
\begin{equation}\label{eqn:newQUE}
 \P\left( \Big | \sum_{i\in I} |u_k(\alpha)| ^2  - \frac{|I|}{N} \Big | \ge \delta\right)\le 
    N^{- \e}/\delta^2.  
\end{equation}
The main problem with this approach is that the obtained QUE is  weak: one would like to replace the above $\e$ with any large $D>0$, as for the GOE in (\ref{eqn:localque}).
For this, it was shown in \cite{BouYauYin2018} that much more general observables than (\ref{eqn:obs1}) also satisfy the eigenvector moment flow parabolic equation (\ref{eqn:EMF}).

These new tractable observables are described as follows. 
Let $I\subset\llbracket 1,N\rrbracket$ be given,  $(\bq_\al)_{\al\in I}$ be any family of fixed vectors, and $C_0\in\mathbb{R}$. Define
\begin{align*}
p_{ij}&=\sum_{\alpha\in I}\langle u_i,\bq_\alpha\rangle\langle u_j,\bq_\alpha\rangle\ \ i\neq j\in\llbracket 1,N\rrbracket,\\
p_{ii}&=\sum_{\alpha\in I}\langle u_i,\bq_\alpha\rangle^2-C_0,\ \ i\in\llbracket 1,N\rrbracket,
\end{align*}
When the $\bq_\alpha$'s are elements of the canonical basis  and $C_0=|I|/N$, this reduces to 
$$
p_{ij}=\sum_{\alpha\in I}u_i(\alpha)u_j(\alpha),\ \ (i\neq j)\,\ 
p_{ii}=\sum_{\alpha\in I}u_i(\alpha)^2-\frac{|I|}{N},\ \ i\in\llbracket 1,N\rrbracket,
$$
and therefore the $p_{ij}$'s become natural partial overlaps measuring quantum unique ergodicity.

\vspace{-0.4cm}
\begin{figure}[h]
\centering
\begin{subfigure}{.4\textwidth}
\centering
\vspace{1.3cm} \hspace{0cm}
\begin{tikzpicture}[scale=0.38]

\draw[fill,black] (1,1) circle [radius=0.2];
\draw[fill,black] (1,2) circle [radius=0.2];
\draw[fill,black] (4,1) circle [radius=0.2];
\draw[fill,black] (4,2) circle [radius=0.2];
\draw[fill,black] (4,3) circle [radius=0.2];
\draw[fill,black] (6,1) circle [radius=0.2];
\draw [thick,->,black] (-2,1) -- (9,1);
\draw [thick,black] (-1,0.8) -- (-1,1.2);
\draw [thick,black] (8,0.8) -- (8,1.2);

\node at (-1,0.4) {\nc 1};
\node at (1,0.4) {\nc$i_1$};
\node at (4,0.4) {\nc $i_2$};
\node at (6,0.4) {\nc $i_3$};
\node at (8,0.4) {\nc $n$};
  
\end{tikzpicture}
\vspace{0.1cm}
  \caption{A configuration $\boeta$ with $\mathcal{N}(\boeta)=6$, $\eta_{i_1}=2$, $\eta_{i_2}=3$, $\eta_{i_3}=1$.
 }
   \label{fig:sub1}
\end{subfigure}
\hspace{1cm}
\begin{subfigure}{.4\textwidth}
  \centering
\begin{tikzpicture}[scale=0.38]
  
\draw[fill,black] (1,1) circle [radius=0.2];
\draw[fill,black] (1,2) circle [radius=0.2];
\draw[fill,black] (1,3) circle [radius=0.2];
\draw[fill,black] (1,4) circle [radius=0.2];
\draw[fill,black] (4,1) circle [radius=0.2];
\draw[fill,black] (4,2) circle [radius=0.2];
\draw[fill,black] (4,3) circle [radius=0.2];
\draw[fill,black] (4,4) circle [radius=0.2];
\draw[fill,black] (4,5) circle [radius=0.2];
\draw[fill,black] (4,6) circle [radius=0.2];
\draw[fill,black] (6,1) circle [radius=0.2];
\draw[fill,black] (6,2) circle [radius=0.2];
\draw [thick,->,black] (-2,1) -- (9,1);
\draw [thick,black] (-1,0.8) -- (-1,1.2);
\draw [thick,black] (8,0.8) -- (8,1.2);

\node at (-1,0.4) {\nc 1};
\node at (1,0.4) {\nc $i_1$};
\node at (4,0.15) {\nc $i_2$};
\node at (6,0.4) {\nc $i_3$};
\node at (8,0.4) {\nc $n$};

\draw [ultra thick,black] (1,1) to[out=180,in=180] (1,3);
\draw [ultra thick,black] (1,2) to[out=-70,in=-90] (6,2);
\draw [ultra thick,black] (1,4) to[out=-70,in=120] (4,1);
\draw [ultra thick,black] (4,6) to[out=-10,in=0] (6,1);
\draw [ultra thick,black] (4,3) to[out=-10,in=0] (4,4);
\draw [ultra thick,black] (4,2) to[out=-10,in=0] (4,5);

\end{tikzpicture}
  \caption{A perfect matching $G\in\mathcal{G}_{\boeta}$. Here, $P(G)=p_{i_1i_1}p_{i_1i_2}p_{i_2i_2}^2p_{i_2i_3}p_{i_3i_1}$.}
  \label{fig:sub1}
\end{subfigure}%

\label{Fig1}
\end{figure}

For any given configuration $\boeta$ as given before (\ref{eqn:obs1}), consider the set of vertices
$
\mathcal{V}_{\boeta}=\{(i,a): 1\leq i\leq n, 1\leq a\leq 2\eta_i\}.
$
Let $\mathcal{G}_{\boeta}$ be the set of perfect matchings of the complete graph on $\mathcal{V}_{\boeta}$, i.e. this is the set of graphs $G$ with vertices $V_{\boeta}$ and edges $\mathcal{E}(G)\subset\{\{v_1,v_2\}: v_1\in \mathcal{V}_{\boeta},v_2\in\mathcal{V}_{\boeta},v_1\neq v_2\}$ being a partition of $\mathcal{V}_{\boeta}$.
For any given edge $e=\{(i_1,a_1),(i_2,a_2)\}$, we define  
$p(e)=p_{i_1,i_2}$, 
$
P(G)=\prod_{e\in \mathcal{E}(G)}p(e)
$ and
\begin{equation}\label{feq}
\widetilde f_{\bla, t}(\boeta)=\frac{1}{\mathcal{M}(\boeta)}\E\left(\sum_{G\in\mathcal{G}_{\boeta}} P(G)\mid \bla\right),\ \ \mathcal{M}(\boeta)=\prod_{i=1}^n (2\eta_i)!!,
\end{equation}
where $(2m)!!=\prod_{k\leq 2m,k \, {\rm odd}}k$.
The following lemma is a key combinatorial fact.

\begin{lemma}
The above function $\widetilde f$ satisfies the eigenvector moment flow equation (\ref{eqn:EMF}).
\end{lemma}

\noindent This new class of observables (\ref{feq})  widely generalizes (\ref{eqn:obs1}) and directly encodes the ${\rm L}^2$ mass of eigenvectors, contrary to (\ref{eqn:obs1}). Together with the above lemma, one can derive a new strong estimate: for a wide class of mean-field models, (\ref{eqn:newQUE}) now holds for arbitrarily {\it large} $\e$.
The mean-field reduction strategy can now be applied in an  efficient way: union bounds are costless thanks to the new small error term.\\

For $d=2,3$, the described mean-field reduction together with the strong version of the eigenvector moment flow should apply to give delocalization in some polynomial regime $W\gg N^{1-\e}$ for some explcit $\e>0$. 
However, this is far from the conjectures from (\ref{eqn:transition exponents}). To approach these transitions, one needs to take into account the geometry of $\mathbb{Z}^d$.

\subsection{Quantum unique ergodicity and the Gaussian free field.}\ At the heuristic level, the QUE method suggests the transition values $W_c$ from (\ref{eqn:transition exponents}). 
More precisely, 
consider a given eigenvector $u=u_k$ associated to a bulk eigenvalue $\lambda_k$. For notational convenience, assume the model's band width is $2W$ instead of $W$.

For $\mathds{1}=(1,\dots,1)\in\mathbb{Z}^d$, define $\mathscr{W}=\llbracket 1,N\rrbracket^d\cap(2W\mathbb{Z}^d+W\mathds{1})$. For any $w\in\mathscr{W}$, let $\mathscr{C}_w=\{\alpha\in\mathbb{Z}^d:
\|w-\alpha\|_\infty\leq W\}$ be the cell of side length $2W$ around $w$.

Let $X_w=\sum_{\alpha\in \mathscr{C}_w} u(\alpha)^2$. Consider a set $I$, $|I|=2^d$, such that the cells $(\mathscr{C}_{w})_{w\in I}$ form a cube $\mathscr{H}$ of size $(4W)^d$.
 Assume one can apply the strong QUE statement (\ref{eqn:localque}) to a Schur complement $Q_e$ of type (\ref{eqn:complement}) where
 $A$ is now chosen to be the $(4W)^d\times (4W)^d$ mean-field matrix indexed by the vertices from $\mathscr{H}$. We would obtain, for any two adjacent cells $\mathscr{C}_{w}, \mathscr{C}_{v}$ with $w,v\in I$,
\begin{equation}\label{eqn:error}
\sum_{\alpha\in\mathscr{C}_{w}}u(\alpha)^2=\sum_{\alpha\in\mathscr{C}_{v}}u(\alpha)^2+\OO\left(N^\e\frac{W^{d/2}}{N^d}\right)
\end{equation}
with overwhelming probability. 
By patching these estimates over successive adjacent cells, this gives
$$
\sum_{\alpha\in\mathscr{C}_{w}}u(\alpha)^2=\left(\frac{W}{N}\right)^d+\OO\left(N^\e\frac{W^{d/2}}{N^d}\right)\times \left(\frac{N}{W}\right),
$$
because there is a path of length $\OO\left(\frac{N}{W}\right)$ between any two cells.
The leading order of $\sum_{\alpha\in\mathscr{C}_{w}}u(\alpha)^2$ is identified (i.e. QUE holds) for $W\gg N^{\frac{2}{d+2}}$. This criterion, improving with the dimension $d$, is more restrictive than (\ref{eqn:transition exponents}) and omits the important fact that the error term in (\ref{eqn:error}) has a random sign.

One may assume that such error terms are asymptotically jointly Gaussian and independent for different pairs of adjacent cells (or at least for sufficiently distant cells). We consider the graph with vertices $\mathscr{W}$ and edges the set of pairs $(v,w)$
such that  $\mathscr{C}_v$ and $\mathscr{C}_w$ are adjacent cells. 
A good model for $(X_w)_{w\in\mathscr{W}}$ therefore is a Gaussian vector such that the increments $X_v-X_w$ are independent, with distribution $\mathscr{N}(0,W^{d}/N^{2d})$ when $(v,w)$ is an edge, and conditioned to (1)
$\sum(X_{v_{i+1}}-X_{v_{i}})=0$ for any closed path $v_1,v_2,\dots, v_j,v_1$ in the graph, (2) $X_{v_0}=(W/N)^d$ to fix the ambiguity about definition of $X$ modulo a constant.
This model is simply the Gaussian free field, with density for $(X_v)_v$ proportional to
$$
e^{-\frac{N^{2d}}{2 W^{d}}\sum_{v\sim w} (x_v-x_w)^2}.
$$
As is well known, the Gaussian free field $(Y_v)_v$ on $\llbracket 1,n\rrbracket^d$ with density $e^{-\frac{1}{2}\sum_{v\sim w} (y_v-y_w)^2}$ conditioned to $Y_{v_0}=0$ has the following typical fluctuation scale, for any deterministic $v$ chosen at macroscopic distance from $v_0$ (see e.g. \cite{Bis}):
$$
{\rm Var}(Y_v)^{1/2}\approx\left\{
\begin{array}{ll}
n^{1/2}&\mbox{for}\,\, d=1,\\
(\log n)^{1/2}&\mbox{for}\,\, d=2,\\
\OO(1)&\mbox{for}\,\, d=3.
\end{array}
\right.
$$
We expect that quantum unique ergodicity (and GOE statistics by the mean-field reduction) holds when ${\rm Var}(X_v)^{1/2}\ll \mathbb{E}(X_v)$.
With $n=N/W$, this means $\frac{W^{d/2}}{N^d}{\rm Var}(Y_v)^{1/2}\ll \frac{W^d}{N^d}$, i.e. $W\gg N^{1/2}$ for $d=1$, $(\log N)^{1/2}$ for $d=2$, $\OO(1)$ for $d=3$.\\

{\it Acknowledgement.}\ The author's knowledge of this topic comes from collaborations with Laszlo Erd{\H o}s,  Horng-Tzer Yau, and Jun Yin. This note reports on joint progress with these authors.

\end{document}